\documentclass[a4paper,12pt,reqno]{amsart}
\author{Wolfgang L\"uck}
        \address{Mathematisches Institut der Universit\"at Bonn\\
                Endenicher Allee 60\\
                53115 Bonn, Germany}
         \email{wolfgang.lueck@him.uni-bonn.de}
         \urladdr{http://www.him.uni-bonn.de/lueck}
\thanks{The first author has been funded by the ERC Advanced Grant 
``KL2MG-interactions'' (no. 662400)  granted by the European Research 
Council, and is funded by the Deutsche Forschungsgemeinschaft (DFG, German Research 
Foundation) under Germany's Excellence Strategy --- GZ 2047/1, Projekt-ID 
390685813, Hausdorff Center for Mathematics at Bonn.}

\author{Bob Oliver}
\address{Universit\'e Sorbonne Paris Nord, LAGA, UMR 7539 du CNRS, 
99, Av. J.-B. Cl\'ement, 93430 Villetaneuse, France.}
\email{bobol@math.univ-paris13.fr}
              \urladdr{https://www.math.univ-paris13.fr/~bobol/}
\thanks{The second author is partially supported by UMR 7539 of the CNRS}

\renewcommand{\today}{\the\day/\the\month/\the\year}

\makeatletter
\@namedef{subjclassname@2020}{%
\textup{2020} Mathematics Subject Classification}
\makeatother

\usepackage{amssymb}
\usepackage{mathrsfs}
\usepackage{enumerate}
\usepackage[all]{xy}
\usepackage[usenames,dvipsnames]{color}
\usepackage[colorlinks=true, citecolor=blue, linkcolor=blue]{hyperref}
\usepackage{datetime2}  

\definecolor{darkgreen}{rgb}{0,0.5,0}
\definecolor{bluegreen}{rgb}{0,0.2,0.8}
\definecolor{darkred}{rgb}{0.8,0,0}
\definecolor{newercolor}{rgb}{0.2,0,1}
\definecolor{darkyellow}{rgb}{0.7,0.7,0}
\definecolor{darkorange}{rgb}{0.8,0.4,0}

\renewenvironment{enumerate}[1][]
{\begin{enumerat}[#1]\setlength{\itemsep}{6pt}}{\end{enumerat}}

\newenvironment{enuma}{\begin{enumerate}[{\rm(a) }]}{\end{enumerate}}

\renewenvironment{itemize}
{\begin{itemiz}\setlength{\itemsep}{6pt} \setlength{\itemindent}{-20pt} }
{\end{itemiz}}

\let\Gamma=\varGamma
\let\Delta=\varDelta

\numberwithin{table}{section}

\newcommand{\boldd}[1]{{\mathversion{bold}\textbf{#1}}}

\newlength{\short}
\newcommand{\4}[1]{\widebar{#1}}
\newcommand{\5}[1]{\widehat{#1}}

\def\pair[#1,#2]{[\hskip-1.5pt[#1,#2]\hskip-1.5pt]}

\SelectTips{cm}{10} \UseTips   

\let\oldcirc=\circ
\renewcommand{\circ}{\mathchoice
    {\mathbin{\scriptstyle\oldcirc}}{\mathbin{\scriptstyle\oldcirc}}
    {\mathbin{\scriptscriptstyle\oldcirc}}
    {\mathbin{\scriptscriptstyle\oldcirc}}}

\def\beq#1\eeq{\begin{equation*}#1\end{equation*}}
\def\beqq#1\eeqq{\begin{equation}#1\end{equation}}

\numberwithin{equation}{section}

\newtheorem{Thm}{Theorem}[section]
\newtheorem{Prop}[Thm]{Proposition}
\newtheorem{Cor}[Thm]{Corollary}
\newtheorem{Lem}[Thm]{Lemma}

\newtheorem{Ex}[Thm]{Example}

\theoremstyle{definition}
\newtheorem{Defi}[Thm]{Definition}

\newcommand{\widebar}[1]
      {\overset{{\mskip3mu\leaders\hrule height0.4pt\hfill\mskip3mu}}{#1}
      \vphantom{#1}}

\newcounter{let} 
\loop\stepcounter{let}
\expandafter\edef\csname cal\alph{let}\endcsname%
{\noexpand\mathcal{\Alph{let}}}
\ifnum\thelet<26\repeat

\loop\stepcounter{let}
\expandafter\edef\csname scr\alph{let}\endcsname%
{\noexpand\mathscr{\Alph{let}}}
\ifnum\thelet<26\repeat

\newcommand{\tdef}[2][]{\expandafter\newcommand\csname#2\endcsname%
{#1\textup{#2}}}
\tdef{Iso}   \tdef{Aut}    \tdef{Out}    \tdef{Inn}    \tdef{Hom}
\tdef{End}   \tdef{Inj}    \tdef{map}    \tdef{Ker}    \tdef{Ob}
\tdef{Mor}   \tdef{Res}    \tdef{Id}     \tdef{Fr}     \tdef{Spin} 
\tdef{rk}    \tdef{conj}   \tdef{incl}   \tdef{proj}   \tdef{diag} 
\tdef{trf}   \tdef{Sol}    \tdef{He}     \tdef{Sz}     \tdef{cj}
\tdef{Rep}   \tdef{pr}    \tdef{Inndiag} \tdef{Outdiag}  \tdef{expt}
\tdef{supp}  \tdef{Isom}   \tdef{ord}    \tdef{Coker}   \tdef{Tr}
\tdef[_]{typ} \tdef[^]{op} \tdef[^]{ab}   \tdef{lcm}  \tdef{McL}
\tdef{restr}  \tdef{Comp}   \tdef{HS}    \tdef{Wh}   \tdef{Gal} 
\tdef{Irr}    \tdef{nr}    \tdef{Suz}    \tdef{Ly}    \tdef{Ru} 
\tdef{Cl}    \tdef{Fix}

\newcommand{\ON}{\textup{O'N}}

\newcommand{\fdef}[1]{\expandafter\newcommand\csname#1\endcsname%
{\mathfrak{#1}}}
\fdef{X}  \fdef{red}  \fdef{foc}  \fdef{hyp}  \fdef{Lie} \fdef{Y}

\newcommand{\bbdef}[1]{\expandafter\newcommand%
\csname#1\endcsname{\mathbb{#1}}}
\bbdef{C} \bbdef{F} \bbdef{R} \bbdef{Z} \bbdef{N} \bbdef{Q} \bbdef{K}

\newcommand{\itdef}[1]{\expandafter\newcommand\csname#1\endcsname%
{\textit{#1}}}
\itdef{PSL}  \itdef{PSU}  \itdef{SL}  \itdef{SU}  \itdef{GL} \itdef{GU}
\itdef{Sp}   \itdef{PSp} \itdef{PSO} \itdef{SO}   \itdef{SD} \itdef{PGU} 
\itdef{PGL}  \itdef{Co}  \itdef{Fi}  \itdef{GO}   \itdef{UT}

\renewcommand{\*}{\mathop{\lower6pt\hbox{\Large{\textup{*}}}}\nolimits}

\newcommand{\sminus}{\smallsetminus}
\newcommand{\lie}[3]{\def\test{#2}\def\tst{G}\ifx\test\tst{{}^{#1}#2_{#3}}
\else{{}^{#1}\!#2_{#3}}\fi}
\newcommand{\syl}[2]{\textup{Syl}_{#1}(#2)}

\newcommand{\defeq}{\overset{\textup{def}}{=}}

\newcommand{\mxfoura}[8]{\left(\begin{smallmatrix}#1&#2&#3&#4\\#5&#6&#7&#8}
\newcommand{\mxfourb}[8]{\\#1&#2&#3&#4\\#5&#6&#7&#8\end{smallmatrix}\right)}

\renewcommand{\:}{\colon}

\newcommand{\nsg}{\trianglelefteq}

\newcommand{\til}[1]{\widetilde{#1}}
\let\too=\longrightarrow

\newcommand{\gen}[1]{{\langle}#1{\rangle}}

\newcommand{\longleft}[1]{\;{\leftarrow%
\count255=0 \loop \mathrel{\mkern-6mu}%
    \relbar\advance\count255 by1\ifnum\count255<#1\repeat}\;}
\newcommand{\longright}[1]{\;{\count255=0 \loop \relbar\mathrel{\mkern-6mu}%
    \advance\count255 by1\ifnum\count255<#1\repeat\rightarrow}\;}
\newcommand{\Right}[2]{\overset{#2}{\longright#1}}
\newcommand{\RIGHT}[3]{\mathrel{\mathop{\kern0pt\longright#1}
        \limits^{#2}_{#3}}}

\newcommand{\LEFT}[3]{\mathrel{\mathop{\kern0pt\longleft#1}\limits^{#2}_{#3}}
}
\newcommand{\dRIGHT}[3]{\mathrel{%
   \mathop{\vcenter{\baselineskip=0pt\hbox{$\kern0pt\longright#1$}%
   \hbox{$\kern0pt\longright#1$}}}\limits^{#2}_{#3}}}
\newcommand{\LRIGHT}[3]{\mathrel{%
   \mathop{\vcenter{\baselineskip=0pt\hbox{$\kern0pt\longleft#1$}%
   \hbox{$\kern0pt\longright#1$}}}\limits^{#2}_{#3}}}
\newcommand{\RLEFT}[3]{\mathrel{%
   \mathop{\vcenter{\baselineskip=0pt\hbox{$\kern0pt\longright#1$}%
   \hbox{$\kern0pt\longleft#1$}}}\limits^{#2}_{#3}}}
\newcommand{\onto}[1]{\;{\count255=0 \loop \relbar\mathrel{\mkern-6mu}%
    \advance\count255 by1
    \ifnum\count255<#1 \repeat \twoheadrightarrow}\;}

 \newcommand{\version}[1] 
     {\begin{center} last edited on #1\\
         last compiled on \today \ at \DTMcurrenttime.\\
         name of tex-file: \jobname
       \end{center}}

     \DeclareMathAlphabet{\matheurm}{U}{eur}{m}{n}
     \newcommand{\FIN}{\mathcal{FIN}}
     \newcommand{\calfj}{\mathcal{F}\!\mathcal{J}}
      \newcommand{\IQ}{\mathbb{Q}}
     \newcommand{\IZ}{\mathbb{Z}}
     \newcommand{\Sub}{\matheurm{Sub}}
     \newcommand{\SubGF}[2]{\Sub_{#2}(#1)}
     \newcommand{\EGF}[2]{E_{#2}(#1)}

     \DeclareMathOperator{\colim}{colim}

\begin{document}


\title{Some criteria concerning the rational vanishing of Whitehead groups}

\begin{abstract}
We give several examples of finite groups $G$ for which the rank of the 
tensor product $\mathbb{Z}\otimes_{\mathbb{Z}\textup{Aut}(G)} 
\textup{Wh}(G)$ is or is not zero. This is motivated by an earlier theorem 
of the first author, which implies as a special case that when this group 
has nonzero rank, the Whitehead group of any other group (finite or 
infinite) that contains $G$ as a normal subgroup is rationally nontrivial. 
\end{abstract}

\maketitle



  \section*{Introduction}\label{sec:introduction}

This paper is motivated by the following question: for which groups 
$\Gamma$ does the Whitehead group $\Wh(\Gamma)$ vanish, integrally or 
rationally? When $\Gamma$ is finite, $\Wh(\Gamma)$ is always finitely 
generated, and its rank is determined in a theorem of Bass (see 
\cite[Theorem 5]{Bass(1964)}). So it is natural to begin studying the 
Whitehead group of an infinite group $\Gamma$ (integrally or rationally) by 
trying to compare it to the Whitehead groups of its finite subgroups.

One means of doing this is provided by the Farrell-Jones Conjecture for the 
algebraic $K$-theory of group rings. This conjecture is known to hold for a 
rather large class of groups, but is open in general. Among other things, 
it makes some predictions about the contributions of finite subgroups of a 
group $\Gamma$ to $\IQ \otimes_{\IZ} \Wh(\Gamma)$ (see Theorem 
\ref{the:criteria_for_rational_vanishing_in_low_and_middle-dimension} below 
for one example of such results). Other known results, such as Theorem 
\ref{the:colim_over_finite_subgroups}, depend instead on assumptions about 
the homology of centralizers of finite cyclic subgroups of $\Gamma$. But so 
far, there are very few results about Whitehead groups known to hold in all 
cases. 

One exception to this is Theorem~\ref{the:Wh(H)_to_Wh(G)_H_normal_finite} 
below. This had originally been predicted by the Farrell-Jones Conjecture, 
but is now known to be true for all groups --- independently of whether or 
not the conjecture holds. This theorem in turn implies 
Corollary~\ref{cor:consequence_of_the:Wh(H)_to_Wh(G)_H_normal_finite}%
\eqref{cor:consequence_of_the:Wh(H)_to_Wh(G)_H_normal_finite:aut}, which 
says, for a group $\Gamma$ with $\IQ\otimes_{\IZ}\Wh(\Gamma)=0$, that 
$\Gamma$ can contain a finite group $H$ as a normal subgroup only if $\IZ 
\otimes_{\IZ[\Aut(H)]} \Wh(H)$ is finite. Since the theorem and corollary 
hold for all groups, this could give some evidence that the Farrell-Jones 
Conjecture for the algebraic $K$-theory of group rings holds more 
generally. 

In Section \ref{s:known}, we describe some of this background in more 
detail, ending with the statements of 
Theorem~\ref{the:Wh(H)_to_Wh(G)_H_normal_finite} and Corollary 
\ref{cor:consequence_of_the:Wh(H)_to_Wh(G)_H_normal_finite}. This then 
motivates the results in the next two sections. In Section 
\ref{s:Whitehead}, we give a general formula for the rank of 
$\Z\otimes_{\Z\Aut(G)}\Wh(G)$, for a finite group $G$, in terms of numbers 
of classes of elements of $G$ under certain equivalence relations. This 
formula is then applied in Section \ref{s:examples}, where we give a wide 
range of examples of finite groups $G$ for which 
$\rk(\Z\otimes_{\Z\Aut(G)}\Wh(G))\ne0$, and hence of groups that cannot 
occur as normal subgroups in a group $\Gamma$ with 
$\IQ\otimes_{\IZ}\Wh(\Gamma)=0$.




\section{Some known results}%
\label{s:known}

We begin by describing some general information about the Whitehead group $G$ of a 
(discrete) group $\Gamma$, and the relation between $\Wh(\Gamma)$ and the 
Whitehead groups of finite subgroups of $\Gamma$.

A group $\Gamma$ is called a \emph{Farrell-Jones group} if it satisfies the 
``Full Farrell-Jones Conjecture'', as formulated, for example, in 
\cite[Conjecture 13.30]{Lueck(2022book)}. The full statement is quite 
complicated and involves $L$-groups as well as $K$-theory, but one easily 
stated special case is that the Whitehead group of a torsion free group is 
always trivial.

The class $\calfj$ of Farrell-Jones groups is quite large, and in fact, no 
groups are known \emph{not} to be in the class. It is known to contain all 
hyperbolic groups, finite-dimensional CAT(0)-groups, solvable groups, 
fundamental groups of manifolds of dimension $\le 3$, and any lattice in a 
locally compact second countable Hausdorff group. Also, it is closed under 
taking subgroups, finite free products, finite direct products, and 
colimits of directed systems (with arbitrary structure maps). For more 
information about the Full Farrell-Jones Conjecture and its consequences, 
we refer to \cite[Chapter~15]{Lueck(2022book)}, and for a description of 
what is currently known about the class $\calfj$, to 
\cite[Chapter~16]{Lueck(2022book)}. 

As one simple example of the role played by Farrell-Jones groups, we note 
the following criterion for the rational vanishing of $K$-theory in degree 
at most $1$. When $H\le\Gamma$ are groups, we let $\Aut_\Gamma(H)$ be the 
group of automorphisms of $H$ of the form $(x\mapsto gxg^{-1})$ for $g\in 
N_\Gamma(H)$.

\begin{Thm}[{\cite[Theorem~17.4]{Lueck(2022book)}}]
\label{the:criteria_for_rational_vanishing_in_low_and_middle-dimension}%
Let $\Gamma$ be a Farrell-Jones group. Consider the following conditions on 
$\Gamma$:
\begin{enuma}

\item[\rm(P)] The  order of every finite cyclic subgroup $C\le\Gamma$ is a prime 
power.

\item[\rm(A)] For every finite cyclic subgroup $1\ne C\le \Gamma$, 
the automorphism group $\Aut(C)$ is generated by $\Aut_\Gamma(C)$ and the 
automorphism $(x\mapsto x^{-1})$.
\end{enuma}
Then 
\begin{enuma}

\item\label{the:criteria_for_rational_vanishing_in_low_and_middle-dimension:n_le_-2}
$K_n(\IZ \Gamma) = 0$ for $n \le -2$;

\item\label{the:criteria_for_rational_vanishing_in_low_and_middle-dimension:n_is-1}
$\IQ \otimes_{\IZ} K_{-1}(\IZ \Gamma) = 0$ if and only if condition 
{\rm(P)} holds; 

\item\label{the:criteria_for_rational_vanishing_in_low_and_middle-dimension:n_is_0}
$\IQ \otimes_{\IZ} \widetilde{K}_{0}(\IZ \Gamma) = 0$ if condition 
{\rm(P)}; and 

\item\label{the:rational_vanishing_of_Wh(G)}
conditions {\rm(P)} and {\rm(A)} $\implies$ $\IQ \otimes_{\IZ} \Wh(\Gamma) = 0$ 
$\implies$ condition {\rm(A)}. 


\end{enuma}
\end{Thm}

The next theorem illustrates the sort of homological condition that can be 
used to get such results, without assuming $G$ is in $\calfj$. 

\begin{Thm}[{\cite[Theorem~1.1]{Lueck-Reich-Rognes-Varisco(2017)}}]
\label{the:colim_over_finite_subgroups}
Let $\Gamma$ be a group. Assume, for every finite cyclic subgroup $C\le 
\Gamma$, that the homology groups $H_1(BC_\Gamma(C);\IZ)$ and 
$H_2(BC_\Gamma(C);\IZ)$ of the centralizer $C_\Gamma(C)$ are finitely 
generated.  Then the  canonical map
  \[
    \colim_{H\in\SubGF{\Gamma}{\FIN}}    \IQ\otimes_{\IZ} \Wh(H) 
    \Right4{} \IQ\otimes_{\IZ} \Wh(\Gamma)
  \]
is injective.
\end{Thm}

Here $\SubGF{\Gamma}{\FIN}$ is the category whose objects are the finite 
subgroups of $\Gamma$, and where a morphism from $H$ to $K$ is the class, 
modulo $\Inn(K)$, of a group homomorphism $f \: H \too K$ of the form 
$(x\mapsto gxg^{-1})$ for some $g\in\Gamma$.

Note that the homology condition appearing in 
Theorem~\ref{the:colim_over_finite_subgroups} holds whenever there is a 
model for the classifying space $\EGF{\Gamma}{\FIN}$ for proper 
$\Gamma$-actions whose $2$-skeleton is cocompact (i.e., the orbit space of 
the $2$-skeleton is finite). Namely, given a finite subgroup $H \subseteq 
G$, we conclude from~\cite[Lemma~1.3]{Lueck(2000a)}, that for every finite 
subgroup $H \subseteq G$ there is a model for $BW_{\Gamma}H$ with finite 
$2$-skeleton.  Since $C_{\Gamma}H$ has finite index in $N_{\Gamma}H$ and 
$W_{\Gamma}H := N_{\Gamma}H/H$ is a finite quotient of $N_{\Gamma}H$, there 
is a model for $BC_{\Gamma}H$ with finite $2$-skeleton.  Examples of groups 
which satisfy the hypotheses of Theorem 
\ref{the:colim_over_finite_subgroups} but are not known to be Farrell-Jones 
groups are $\Out(F_n)$ (where $F_n$ is a free group on $n$ letters) and 
Thompson's group.

In contrast to Theorems 
\ref{the:criteria_for_rational_vanishing_in_low_and_middle-dimension} and 
\ref{the:colim_over_finite_subgroups}, the next theorem holds for 
\emph{all} (discrete) groups. Let $i\:G \to \Gamma$ be the inclusion 
of a finite normal subgroup $G$ in a group $\Gamma$, and let $i_*\: \Wh(G) \too 
\Wh(\Gamma)$ be the induced homomorphism.  The conjugation actions of 
$\Gamma$ on $G$ and $\Gamma$ induce $\Gamma$-actions on 
$\Wh(G)$ and $\Wh(\Gamma)$, of which the latter is trivial. Hence $i_*$ induces a 
homomorphism 
	\[ \overline{\imath_*} \colon \IZ \otimes_{\IZ \Gamma} \Wh(G) 
	\Right4{} \Wh(\Gamma). \]

\begin{Thm}[{\cite[Theorem~9.38]{Lueck(2002)}}]
\label{the:Wh(H)_to_Wh(G)_H_normal_finite}
Let $i\colon  G \to \Gamma$ be the inclusion of a finite normal subgroup 
$G$ into an arbitrary group $\Gamma$. Then the homomorphism 
	\[ \overline{\imath_*} \: \Z\otimes_{\Z\Gamma}\Wh(G) \Right4{} 
	\Wh(\Gamma) \]
induced by $i$ has finite kernel.
\end{Thm}

\begin{Cor}\label{cor:consequence_of_the:Wh(H)_to_Wh(G)_H_normal_finite}
Let $\Gamma$ be a group such that $\IQ \otimes_{\IZ} \Wh(\Gamma)=0$. 
Then the following hold.
\begin{enuma}
    
\item\label{cor:consequence_of_the:Wh(H)_to_Wh(G)_H_normal_finite:center} 
Every finite cyclic subgroup of the center $Z(\Gamma)$ has order $1$, $2$, $3$, 
$4$, or $6$.

\item\label{cor:consequence_of_the:Wh(H)_to_Wh(G)_H_normal_finite:aut} For 
each finite normal subgroup $G\nsg \Gamma$, the group $\IZ \otimes_{\IZ 
\Gamma} \Wh(G)$ and hence also $\IZ \otimes_{\IZ\Aut(G)} \Wh(G)$ is 
finite, where $\Aut(G)$ acts in the canonical way on $\Wh(G)$ and 
trivially on $\IZ$.

\end{enuma}
\end{Cor}

\begin{proof} By a theorem of Bass \cite[\S\,7]{Bass(1964)}, the Whitehead 
group of a finite group is always finitely generated. So (b) follows 
immediately from Theorem \ref{the:Wh(H)_to_Wh(G)_H_normal_finite}. 

If $C\le Z(\Gamma)$ is cyclic of finite order, then by Theorem 
\ref{the:Wh(H)_to_Wh(G)_H_normal_finite} again, $\Q\otimes_\Z\Wh(C)=0$, and 
$\Wh(G)$ has rank $0$. So 
by Bass's theorem again (see also Theorem \ref{t:rk(Wh^A)} below), the 
group $C$ has at most two generators, and hence $|C|\le4$ or $|C|=6$.
\end{proof}

Motivated by 
Corollary~\ref{cor:consequence_of_the:Wh(H)_to_Wh(G)_H_normal_finite}%
\eqref{cor:consequence_of_the:Wh(H)_to_Wh(G)_H_normal_finite:aut}, in 
the rest of the paper, we 
look for finite groups $G$ for which $\rk(\IZ \otimes_{\IZ[\Aut(G)]} 
\Wh(G))>0$, since no such group can occur as a normal subgroup in $\Gamma$ 
if $\IQ \otimes_{\IZ} \Wh(\Gamma)=0$. Note that this condition on $\IZ 
\otimes_{\IZ[\Aut(G)]} \Wh(G)$ is independent of $\Gamma$.



\section{Whitehead groups and conjugacy classes}
\label{s:Whitehead}

From now on, we focus attention on Whitehead groups of finite groups. 
The main result in this section is a generalization of a theorem of 
Bass \cite[\S\,7]{Bass(1964)} describing the rank of $\Wh(G)$ for finite 
$G$.

When $K$ is a field of characteristic $0$ and $n>0$ is an integer, we let 
$\mu_n$ be the group of $n$-th roots of unity in an algebraic closure of 
$K$, and regard the Galois group $\Gal(K(\mu_n)/K)$ as a subgroup of 
$(\Z/n)^\times$. Thus each automorphism $\gamma\in\Gal(K(\mu_n)/K)$ is 
identified with the unique class $a{+}n\Z\in(\Z/n)^\times$ such 
that $\gamma(\zeta)=\zeta^a$ for $\zeta\in\mu_n$. For example, 
	\[ \Gal(\Q(\mu_n)/\Q) =(\Z/n)^\times, \quad 
	\Gal(\R(\mu_n)/\R) =\{\4{\pm1}\} ~\textup{if $n\ge3$,}
	\quad\textup{and}\quad \Gal(\C(\mu_n)/\C)=1. \]

\begin{Defi}\label{d:K-A-conj}
Let $K$ be a field of characteristic $0$. When $G$ is a finite group, two 
elements $g,h\in G$ of the same order $n$ are \emph{$K$-$G$-conjugate} if 
$g$ is $G$-conjugate to $h^a$ for some $a+n\Z\in\Gal(K(\mu_n)/K)$. More 
generally, when $A\le\Aut(G)$ is a subgroup containing $\Inn(G)$, the elements 
$g,h\in G$ of order $n$ are \emph{$K$-$A$-conjugate} if $g$ is in the 
$A$-orbit of $h^a$ for some $a+n\Z\in\Gal(K(\mu_n)/K)$.
\end{Defi}

Thus $g,h\in G$ are $\R$-$A$-conjugate if $g$ is in the $A$-orbit of $h$ or 
$h^{-1}$, while they are $\Q$-$A$-conjugate if the cyclic subgroups 
$\gen{g}$ and $\gen{h}$ are in the same $A$-orbit. In particular, the 
number of $\Q$-$G$-conjugacy classes is equal to the number of conjugacy 
classes of cyclic subgroups of $G$.

As usual, by the \emph{rank} of a finitely generated abelian group $B$, we 
mean the rank of its free part; i.e., the order of its largest 
$\Z$-linearly independent subset. Thus $\rk(B)=\dim(\Q\otimes_{\Z}B)$. 

For a field $K$ of characteristic $0$ and a finite group $G$, we let 
$\Irr_K(G)$ be the set of isomorphism classes of irreducible $KG$-modules, 
regarded as a finite $\Aut(G)$-set. Via character theory, one shows the 
following:

\begin{Prop}\label{p:K-A-conj}
Let $K$ be a field of characteristic $0$. For each finite group $G$, and 
each group of automorphisms $A\le\Aut(G)$ containing $\Inn(G)$, 
	\[ |\Irr_K(G)/A| = \# \bigl\{\textup{$K$-$A$-conjugacy classes of 
	elements of $G$}\bigr\}. \]
\end{Prop}

\begin{proof} Let $R_K(G)$ be the representation ring for 
$KG$-representations: thus a free abelian group with basis $\Irr_K(G)$. Let 
$\Cl_K(G)$ be the space of all maps $G\too K$ that are constant on 
$K$-$G$-conjugacy classes, regarded as a vector space over $K$.
By~\cite[\S12.4, Corollary 2]{Serre(1977)}, the characters of the elements in 
$\Irr_K(G)$ form a $K$-basis for $\Cl_K(G)$. Hence the sums of the characters in 
each $A$-orbit form a basis for $\Fix(A,\Cl_K(G))$, and so 
	\beq |\Irr_K(G)/A| = \dim_K(\Fix(A,\Cl_K(G))) 
	= \# \bigl\{\textup{$K$-$A$-conjugacy classes of elements of 
	$G$}\bigr\}. \qedhere \eeq
\end{proof}

We refer to \cite[Section I.2a]{Oliver(1988)} for the definition of 
reduced norms for finite dimensional semisimple $\Q$-algebras.

\begin{Prop}\label{p:R/A}
Let $G$ be a finite group. Set $K=Z(\Q G)$, a product of fields, and let 
$R\le K$ be its unique maximal order (the product of the rings of integers in 
the factors). The reduced norm induces a homomorphism
	\[ \nr\: K_1(\Z G) \Right4{} R^\times \]
that commutes with the actions of $\Out(G)$ on $K_1(\Z G)$ and on 
$R^\times$, and whose kernel and cokernel are both finite. In particular, 
for each $A\le\Aut(G)$ containing $\Inn(G)$, 
	\[ \rk(\Z\otimes_{\Z A}\Wh(G)) = \rk(\Fix(A,\Wh(G))) 
	= \rk(\Fix(A,K_1(\Z G))) = \rk(\Fix(A,R^\times)). \]
\end{Prop}

\begin{proof} The fact that the reduced norm induces a homomorphism whose 
kernel and cokernel are finite was shown by Swan (see, e.g., 
\cite[Chapter 8]{Swan(1970)} or \cite[Theorem I.2.5(ii)]{Oliver(1988)}). 

Clearly, $\nr$ commutes with the actions of $A$. So 
$\Fix(A,K_1(\Z G))$ and $\Fix(A,R^\times)$ have the same rank.
\end{proof}

We also need Dirichlet's units theorem in the following form.

\begin{Lem}\label{l:Dirichlet}
Let $K\supseteq\Q$ be a finite extension, and let $R\subseteq K$ be the 
ring of integers. Then 
	\[ \rk(R^\times) = \#\bigl\{\textup{field factors in 
	$\R\otimes_{\Q}K$}\bigr\} 
	- 1. \]
\end{Lem}

\begin{proof} The ring $\R\otimes_{\Q}K$ is isomorphic to a product of one 
copy of $\R$ for each embedding of $K$ into $\R$, and one copy of $\C$ for 
each pair of conjugate embeddings of $K$ into $\C$ with image not in $\R$. 
So the statement follows from Dirichlet's theorem in its usual formulation 
(see \cite[Theorem B.6]{Stewart}). 
\end{proof}

\begin{Prop}\label{p:Wh-Irr(G)}
Let $G$ be a finite group. Then for each $A\le\Aut(G)$ containing $\Inn(G)$, 
we have 
	\[ \rk(\Z\otimes_{\Z A}\Wh(G)) = \rk(\Fix(A,\Wh(G))) = 
	|\Irr_{\R}(G)/A| - |\Irr_{\Q}(G)/A|. \]
\end{Prop}

\begin{proof} Set $K=Z(\Q G)$, a product of fields, and let $R\le K$ be its 
unique maximal order (the product of the rings of integers in the factors). 
Then $\rk(\Fix(A,\Wh(G)))=\rk(\Fix(A,R^\times))$ by Proposition~\ref{p:R/A}, 
and it remains to describe the rank of $\Fix(A,R^\times)$ in terms of 
representations of $G$. 

Set $K_0=\Fix(A,K)$ and $R_0=\Fix(A,R)$: the subgroups of elements fixed by 
$A$. By Wedderburn's theorem, there is a natural bijection of $A$-sets 
from $\Irr_{\Q}(G)$ to the set of simple factors in $\Q G$, and hence 
to the set of field factors in $Z(\Q G)=K$. By a similar argument, there is 
a natural bijection of $A$-sets from $\Irr_{\R}(G)$ to the set of simple 
factors in $\R G$, and hence to the set of field factors in $Z(\R 
G)\cong\R\otimes_{\Q}K$. Also, by Galois theory, if a subgroup $A_0\le A$ 
sends a field factor to itself, then the set of elements in that field 
fixed by $A_0$ is a subfield. Thus the number of field factors in 
$K_0=\Fix(A,K)$ is equal to $|\Irr_{\Q}(G)/A|$, and the number of 
field factors in $\R\otimes_{\Q}K_0=\Fix(A,\R\otimes_{\Q}K)$ is equal to 
$|\Irr_{\R}(G)/A|$. 

Clearly, $R_0$ contains the product of the rings of integers in the field 
factors of $K_0$, with equality since for each field factor $L$ of $K_0$, 
the image of $R_0$ under projection to $L$ is a finitely generated subring 
and hence contained in the ring of integers (see, e.g., \cite[Lemma 
2.8]{Stewart}). So 
by Lemma~\ref{l:Dirichlet}, 
	\begin{align*} 
	\rk(\Fix(A,R^\times)) = \rk((R_0)^\times) &= \#\bigl\{\textup{field 
	factors in $\R\otimes_{\Q}K_0$}\bigr\} - \#\bigl\{\textup{field 
	factors in $K_0$}\bigr\} \\
	&= |\Irr_{\R}(G)/A| - |\Irr_{\Q}(G)/A|. \qedhere
	\end{align*}
\end{proof}

The following theorem is due to Bass \cite[Theorem 5]{Bass(1964)} 
when $A=\Inn(G)$. 

\begin{Thm}\label{t:rk(Wh^A)}
For each finite group $G$, and each subgroup $A\le\Aut(G)$ containing 
$\Inn(G)$, 
	\begin{multline*} 
	\dim_{\Q}\bigl((\Q\otimes_{\Z A}\Wh(G))\bigr) = 
	\#\bigl\{ \textup{$\R$-$A$-conjugacy classes of elements in $G$} 
	\bigr\} \\
	- \# \bigl\{ \textup{$\Q$-$A$-conjugacy classes of elements in $G$} 
	\bigr\}. 
	\end{multline*}
\end{Thm}

\begin{proof} This follows directly from Propositions~\ref{p:Wh-Irr(G)} 
and~\ref{p:K-A-conj}. 
\end{proof}




\section{Examples}
\label{s:examples}

\newcommand{\NN}{{\scrn}}

Throughout this section, it will be convenient to define, for each finite 
group $G$,
	\[ \NN_G = 
	\#\bigl\{\textup{$\R$-$\Aut(G)$-conjugacy classes in $G$}\bigr\} - 
	\#\bigl\{\textup{$\Q$-$\Aut(G)$-conjugacy classes in $G$}\bigr\} . 
	\]
Thus by Theorem \ref{t:rk(Wh^A)}, 
	\[ \NN_G = \rk(\Z\otimes_{\Z\Aut(G)}\Wh(G)) = 
	\dim_{\Q}(\Q\otimes_{\Z\Aut(G)}\Wh(G)). \]
So by Corollary 
\ref{cor:consequence_of_the:Wh(H)_to_Wh(G)_H_normal_finite}%
\eqref{cor:consequence_of_the:Wh(H)_to_Wh(G)_H_normal_finite:aut}, if $\NN_G>0$, 
then $\Q\otimes_{\Z}\Wh(\Gamma)\ne0$ for every group $\Gamma$ that contains $G$ 
as a normal subgroup. 

We now construct a wide variety of examples of finite groups $G$ with 
$\NN_G>0$, including some small metacyclic groups, metacyclic $p$-groups, 
and simple groups. For example, we show that the smallest group with 
$\NN_G>0$ is the nonabelian group of order $55$, and that $\NN_G>0$ when 
$G$ is the nonabelian group of order $p^3$ and exponent $p^2$ and $p\ge5$ 
is prime. We also determine exactly which of the sporadic simple groups, 
and which of the linear groups $\PSL_n(q)$, satisfy $\NN_G>0$.

We first need some tools for constructing automorphisms.

\begin{Lem}\label{l:abel.n}
\begin{enuma} 

\item Let $H\nsg G$ be a pair of finite groups such that $H$ is abelian and 
$G/H$ is cyclic. Then for each $x\in G$ such that $G=H\gen{x}$, and each 
$a\in\Z$ prime to $|G|$ such that $a\equiv1$ (mod $|G/H|$), there is 
$\alpha\in\Aut(G)$ such that $\alpha(h)=h^a$ for each $h\in H\cup\{x\}$. 

\item Let $G$ be a group, let $Z\le Z(G)$ be a central subgroup, and let 
$\psi\in\Hom(G,Z)$ be a homomorphism such that $Z\le\Ker(\psi)$. Then there 
is $\alpha\in\Aut(G)$ such that $\alpha(g)=g\psi(g)$ for each $g\in G$.

\end{enuma}
\end{Lem}

\begin{proof} \textbf{(a) } Define $\alpha\:G\too G$ by setting, 
for each $h\in H$ and each $i\in\Z$, $\alpha(hx^i)=h^ax^{ai}$. By 
assumption, every element in $G$ can be written in this form. If 
$hx^i=kx^j$ for $h,k\in H$ and $i,j\in\Z$, then $k^{-1}h=x^{j-i}$, so 
$k^{-a}h^a=(k^{-1}h)^a=x^{aj-ai}$, and hence $h^ax^{ai}=k^ax^{aj}$. So 
$\alpha$ is well defined as a map of sets. Finally, if $h,k\in H$ and 
$i,j\in\Z$ are arbitrary, then 
	\begin{align*} 
	\alpha((hx^i)(kx^j)) &= \alpha(h(x^ikx^{-i})x^{i+j}) = 
	h^a(x^ik^ax^{-i})x^{ai+aj} = h^a(x^{ai}k^ax^{-ai})x^{ai+aj} \\ 
	&= (h^ax^{ai})(k^ax^{aj}) = \alpha(hx^i)\alpha(kx^j), 
	\end{align*}
where the third equality holds since $x^{(a-1)i}\in H$ and hence commutes 
with $k^a$. So $\alpha$ is an automorphism.

\noindent\textbf{(b) } One easily checks that $\alpha$ is a homomorphism 
with inverse $(g\mapsto g\psi(g)^{-1})$.
\end{proof}

As a first, very simple, application of Lemma \ref{l:abel.n}(a), we have:

\begin{Ex}\label{ex:abel.n}
If a finite group $G$ contains a normal abelian subgroup of index at most 
$3$, then $\NN_G=0$.
\end{Ex}

\begin{proof} Assume $H\nsg G$ is abelian of index at most $3$. Fix $g\in 
G$; then either $g\in H$ or $G=H\gen{g}$. By Lemma \ref{l:abel.n}(a), for each 
$a\in\Z$ such that $\gcd(a,|G|)=1$, and such that $a\equiv1$ (mod $3$) if 
$|G/H|=3$, there is $\alpha\in\Aut(G)$ such that $\alpha(g)=g^a$. Thus 
all generators of $\gen{g}$ are $\R$-$\Aut(G)$-conjugate to $g$. Since 
$g\in G$ was arbitrary, Theorem~\ref{t:rk(Wh^A)} now implies
	\beq \NN_G = \#\bigl\{\textup{$\R$-$\Aut(G)$-conjugacy classes} 
	\bigr\} - 
	\#\bigl\{\textup{$\Q$-$\Aut(G)$-conjugacy classes} \bigr\} = 0. 
	\qedhere \eeq
\end{proof}

With a little more work, one can show that if $G$ is a finite group with 
$\NN_G=0$, then $\NN_{G\times H}=0$ for each finite group $H$ that contains 
an abelian subgroup of index at most $2$. However, if we let $G\cong 
C_5\rtimes C_4$ (induced by an injection $C_4\too\Aut(C_5)$), and let $H$ 
be a nonabelian group of order $21$, then $\NN_G=\NN_H=0$, but 
$\NN_{G\times H}>0$.

We next look at some more small groups $G$ for which $\NN_G$ vanishes.

\begin{Ex} \label{ex:|G|<55}
If $G$ is a group of order at most $54$, then $\NN_G=0$.
\end{Ex}

\begin{proof} If each element of $G$ has order dividing $4$ or $6$, then 
all generators of each cyclic subgroup of $G$ are $\R$-$G$-conjugate, so 
$\NN_G=0$. If $G$ is abelian, or contains a 
normal abelian subgroup of index $2$ or $3$, then $\NN_G=0$ by Example 
\ref{ex:abel.n}. From now on, we use these without repeating the 
references each time. 

When $q=p^k$ for a prime $p$ and $k\ge1$, we let $E_q$ denote an elementary 
abelian $p$-group of order $q$ and rank $k$. When $P$ is a $p$-group for a 
prime $p$, we let $\Phi(P)$ denote its Frattini subgroup: the intersection 
of the maximal proper subgroups of $P$, and the smallest normal subgroup 
such that $P/\Phi(P)$ is elementary abelian.

\noindent\boldd{Case 1: $|G|=n$ where $n\le53$ is odd.} If $n$ is prime or the 
square of a prime, or $n=35$ or $45$, then $G$ is abelian, so $\NN_G=0$. 
Otherwise, $n=3m$ where $m>3$ is prime or $m=9$, so $G$ contains a normal 
abelian subgroup of index $3$, and $\NN_G=0$. 

\noindent\boldd{Case 2: $|G|=2n$ where $n\le27$ is odd.} In these cases, 
$G$ always contains a normal subgroup $H\nsg G$ of order $n$ and index $2$. 
If $H$ is abelian, then $\NN_G=0$. If $H$ is nonabelian, then $n=21$ or 
$27$.

If $n=21$ and $H$ is nonabelian, then $G$ is a semidirect product of the 
form $C_{14}\rtimes C_3$ or $C_7\rtimes C_6$. In the first case, $\NN_G=0$. 
In the second case, all elements of $G$ have order $7$ or a divisor of $6$, 
and the elements of order $7$ are permuted transitively by $\Aut(G)$ by 
Lemma \ref{l:abel.n}(a). So $\NN_G=0$ also in this case.

This leaves the case where $|H|=n=27$ and $H$ is nonabelian. If $H$ has exponent 
$3$, then the order of each element of $G$ divides $6$, and so $\NN_G=0$. 
So assume $H$ is nonabelian of exponent $9$. There are three subgroups of 
order $9$ in $H$, at least one of which must be invariant under the 
conjugation action of $G/H\cong C_2$. So there is a normal subgroup $K\nsg 
G$ with $K\cong C_9$. If $C_G(K)>K$, then $G$ has a normal abelian subgroup 
of index $3$, and so $\NN_G=0$. Otherwise, $G/K\cong C_6$ (since 
$\Aut(K)\cong C_6$), and $G$ is a semidirect product $K\rtimes C_6$ where 
$C_G(K)=K$. 

In this last case, $G$ has presentation 
$\gen{a,b\,|\,a^9=1=b^6,~bab^{-1}=a^2}$. So by Lemma \ref{l:abel.n}(a), 
there is $\alpha\in\Aut(G)$ such that $\alpha(a)=a^7$ and $\alpha(b)=b$. 
Then $\alpha$ induces the identity on $G/\gen{a^3}$, so for each cyclic 
subgroup $C\le G$ of order $9$, $\alpha(C)=C$ and hence all generators of 
$C$ are $\R$-$\Aut(G)$-conjugate. Since all other elements have order 
dividing $6$, this shows that $\NN_G=0$.

\noindent\boldd{Case 3: $|G|=4n$ where $n\le13$ is odd.} If $n=1$ or 
$n\equiv3$ (mod $4$), then $G$ contains an abelian subgroup of index $2$. 
If $n=5$ or $13$, then either $G$ has an abelian subgroup of index $2$, or 
$G\cong C_n\rtimes C_4$, and each element of $G$ has order $1$, $2$, $4$, 
or $n$. All elements in $G$ of order $n$ are $\R$-$\Aut(G)$-conjugate by 
Lemma \ref{l:abel.n}(a), so $\NN_G=0$. 

Assume $n=9$, so $|G|=36$. If there is $H\nsg G$ of order $9$, then 
either $G$ has an abelian subgroup of index $2$, or $H\cong E_9$, 
$C_G(H)=H$, and all elements of $G$ have order $1,2,3,4,6$. So $\NN_G=0$ in 
all such cases. Otherwise, $|\syl3G|=4$, the conjugation action on this set 
induces a homomorphism $\chi\:G\too\Sigma_4$ with image $A_4$, and so $G$ 
has a normal abelian subgroup of index $3$.

\noindent\boldd{Case 4: $|G|=8$, $16$, or $32$.} If $G/Z(G)$ has an element 
of order $8$, generated by the class of $g\in G$, then $Z(G)\gen{g}$ is 
abelian of index at most $2$ in $G$. Also, if $G$ has an element of order 
$16$, then it generates a cyclic subgroup of index at most $2$. So 
$\NN_G=0$ in all such cases. 

Assume from now on that $G/Z(G)$ has exponent at most $4$ and $G$ has 
exponent at most $8$. So if $g\in G$ has order $8$, then $g^4\in Z(G)$. If 
$g\notin\Phi(G)$, then by Lemma \ref{l:abel.n}(b), applied with 
$Z=\gen{g^4}$, there is $\alpha\in\Aut(G)$ such that $\alpha(g)=g^5$. So 
all generators of $\gen{g}$ are $\R$-$\Aut(G)$-conjugate in this case. 

If $g\in\Phi(G)$ has order $8$, then $\Phi(G)=\gen{g}$ and $G/\Phi(G)\cong 
C_2\times C_2$, which is impossible unless $G$ has a cyclic subgroup of 
index $2$.

\noindent\boldd{Case 5: $|G|=24$.} If $|\syl3G|=1$, then $G$ contains a 
normal abelian subgroup of index $2$, so $\NN_G=0$. Otherwise, 
$|\syl3G|=4$, and the conjugation action on this set defines a homomorphism 
$\chi\:G\too\Sigma_4$ whose image contains $A_4$. So either 
$G\cong\Sigma_4$, or $G/Z(G)\cong A_4$ and $|Z(G)|=2$, and in all such 
cases, all elements of $G$ have order dividing $4$ or $6$. So $\NN_G=0$.

\noindent\boldd{Case 6: $|G|=40$.} In all cases, $|\syl5G|=1$. Then either 
$G$ has an abelian subgroup of index $2$ (and $\NN_G=0$), or there is a 
normal subgroup $H\nsg G$ such that $H\cong C_{10}$, $G/H\cong C_4$, and 
$C_G(H)=H$. 

Assume we are in this last case. By Lemma \ref{l:abel.n}(a), the generators 
of $H$ are permuted transitively by $\Aut(G)$, and similarly 
for the elements of order $5$. All other elements have order dividing $8$. 
Let $x\in H$ be the element of order $2$, and let $\alpha\in\Aut(G)$ be the 
automorphism that is the identity on the subgroup of index $2$ that 
contains $H$, and sends $g$ to $gx$ otherwise. If $g\in G$ has order $8$, 
then $x=g^4$ and $\alpha(g)=gx=g^5$, so all generators of $\gen{g}$ are 
$\R$-$\Aut(G)$-conjugate. Thus $\NN_G=0$. 

\noindent\boldd{Case 7: $|G|=48$.} By the Sylow theorems, $|\syl3G|=1$, 
$4$, or $16$. 

\noindent\boldd{Case 7A: } If $|\syl3G|=16$, then the Sylow $2$-subgroup 
$S\le G$ is normal and contains all elements not of order $3$. So $S$ has 
an automorphism of order $3$ that fixes only the identity element, hence 
$[S,S]$ is trivial or noncyclic, which implies $S\cong E_{16}$ or 
$C_4\times C_4$. Thus all elements of $G$ have order at most $4$, 
and so $\NN_G=0$.

\noindent\boldd{Case 7B: } Assume $|\syl3G|=4$, and let 
$\chi\:G\too\Sigma_4$ be the homomorphism defined by the conjugation 
action. If $\chi$ is onto, then $\Ker(\chi)=Z(G)=\gen{x}$ where $|x|=2$, so 
by Lemma \ref{l:abel.n}(b), there is $\alpha\in\Aut(G)$ such that 
$\alpha|_{\chi^{-1}(A_4)}=\Id$ and $\alpha(g)=gx$ for $g\in G$ with 
$\chi(g)\notin A_4$. So if $g\in G$ has order $8$, then $g^4=x$ and 
$\alpha(g)=gx=g^5$. Thus all generators of each cyclic subgroup of order 
$8$ are $\R$-$\Aut(G)$-conjugate. Since all elements of $G$ have order 
dividing $6$ or $8$, this proves that $\NN_G=0$. 

If $\chi(G)=A_4$, then either $\Ker(\chi)\cong E_4$, or 
$\Ker(\chi)=Z(G)\cong C_4$ and $G\cong C_4\times A_4$ or $C_4\circ\SL_2(3)$ 
(central product). In all such cases, all elements of $G$ have order 
dividing $4$ or $6$, so $\NN_G=0$.

\noindent\boldd{Case 7C: } Now assume $|\syl3G|=1$, let $H\nsg G$ be the 
normal subgroup of order $3$, and choose $S\in\syl2G$. Then $C_G(H)=HT$ for 
some $T\le S$ of index at most $2$, and $\NN_G=0$ if $T$ is abelian. So 
assume $S\cong G/H$ is nonabelian of order $16$, with a nonabelian 
subgroup $T\le S$ of index at most $2$. 

If $T=S$, then $G\cong H\times S$. If $S/Z(S)$ has a cyclic subgroup of 
index $2$, then $S$ has an abelian subgroup of index $2$. Otherwise, 
$|Z(S)|=2$ and $S/Z(S)\cong E_8$, in which case some pair of commuting 
elements in $S/Z(S)$ lift to commuting elements in $S$. So in all cases, 
$S$, and hence $G$, have abelian subgroups of index $2$, and so $\NN_G=0$.

Assume for the rest of the proof that $T$ has index $2$ in $S$. Since $G$ 
is a semidirect product $G=H\rtimes S$, each $\alpha\in\Aut(S)$ with 
$\alpha(T)=T$ extends to $\5\alpha\in\Aut(G)$ such that 
$\5\alpha|_H=\Id_H$. Also, there is $\beta\in\Aut(G)$ such that 
$\beta|_S=\Id_S$ and $\beta(h)=h^{-1}$ for $h\in H$.

If $S$ has an element $g$ of order $8$, then since $T$ is nonabelian, we 
have $S=\gen{g,x}$ and $T=\gen{g^2,x}$ for some $x\in T\sminus\gen{g}$. By 
Lemma \ref{l:abel.n}(a), for each odd $a\in\Z$, there is 
$\alpha\in\Aut(S)$ that $\alpha(g)=g^a$ and $\alpha(x)=x^a$, so 
$\alpha(T)=T$, and hence $\alpha$ extends to $\5\alpha\in\Aut(G)$. Thus 
$N_{\Aut(G)}(\gen{g})$ permutes transitively the generators of $\gen{g}$. 

If $g\in G$ has order $12$, then $H=\gen{g^4}$ and $g^3\in T$. Thus 
$\beta(g^4)=g^{-4}$ and $\beta(g^3)=g^3$, so $\beta(g)=g^5$, proving that all 
generators of $\gen{g}$ are $\R$-$\Aut(G)$-conjugate. All elements of $G$ 
have order $12$, $8$, $6$, or at most $4$, so $\NN_G=0$ in all of these 
cases. 
\end{proof}

By Example \ref{ex:|G|<55}, the smallest example of a finite group $G$ such 
that $\NN_G>0$ has order at least $55$. We now prove that there is such a 
group. Let $\varphi$ denote the Euler function: 
$\varphi(n)=|(\Z/n)^\times|$.

\begin{Ex} \label{ex:pm}
Fix a prime $p$, and an integer $m\ge3$ such that $m|(p-1)$. Let $G$ be a 
semidirect product $G\cong C_p\rtimes C_m$, where $C_m$ acts faithfully on 
$C_p$. Then $\NN_G\ge(\varphi(m)/2)-1$, with equality if $m$ is prime. Thus 
$\NN_G>0$ if $m=5$ or $m\ge7$.
\end{Ex}

\begin{proof} Fix elements $a,b\in G$ with $|a|=p$ and $|b|=m$, and set 
$H=\gen{a}\nsg G$. Let $k\in\Z$ be such that $bab^{-1}=a^k$. For each 
$\alpha\in\Aut(G)$, $\alpha(a)=a^s$ for some $s$, and so 
$\alpha(b)\alpha(a)\alpha(b)^{-1}=\alpha(a)^s=b\alpha(a)b^{-1}$. Hence 
$\alpha(b)\in bH$. In other words, each automorphism of $G$ induces the 
identity on $G/H\cong C_m$. 

Thus the $\varphi(m)$ generators of $\gen{b}$ all lie in one 
$\Q$-$\Aut(G)$-conjugacy class, but in separate $\Aut(G)$-conjugacy 
classes. Since each $\R$-$\Aut(G)$-conjugacy class is the union of at most 
two $\Aut(G)$-conjugacy classes, there are at least $\varphi(m)/2$ 
$\R$-$\Aut(G)$-conjugacy classes in the $\Q$-$\Aut(G)$-conjugacy class of 
$b$. So $\NN_G\ge(\varphi(m)/2)-1$. 

If $m$ is prime, then every nonidentity 
element of $G$ has order $p$ or $m$ and hence is $\Q$-$G$-conjugate to $a$ 
or $b$. All generators of $\gen{a}$ are $\Aut(G)$-conjugate to $a$, and there 
are $\varphi(m)/2$ $\R$-$\Aut(G)$-conjugacy classes of generators 
of $\gen{b}$. So $\NN_G=(\varphi(m)/2)-1$ in this case.
\end{proof}

With a little more work, one can show that 
$\NN_G=\sum_{2<d|m}((\varphi(d)/2)-1)$ in the situation of Example 
\ref{ex:pm}. 

The next example also involves metacyclic groups, and shows that there are 
$p$-groups $G$ (for an arbitrary prime $p$) such that $\NN_G>0$. Note that 
the smallest examples constructed in this way have order $p^3$ when 
$p\ge5$, or order $3^6$ or $2^{10}$ when $p=3$ or $2$.

\begin{Ex} \label{ex:p-gp}
Fix integers $1<r|q$, and let $G$ be the group of order $qr^2$ with 
presentation
	\[ G = \gen{a,b \,|\, a^{qr}=1=b^r,~ bab^{-1}=a^{q+1} }. \]
Then $\NN_G>0$ if $\varphi(r)>4$, or if $\varphi(r)=4$ and $q$ is odd or 
$2r|q$. 
\end{Ex}

\begin{proof} Set $t=q+1$ for short. Note first that 
$t^s=(1+q)^s\equiv1+sq$ (mod $q^2$) for each $s\ge1$. In particular, 
$t^r\equiv1$ (mod $qr$), so the above presentation does define a group of 
order $qr^2$. 

Let $\alpha\in\Aut(G)$ be an automorphism that normalizes the cyclic 
subgroup $\gen{b}$. Thus 
	\[ \alpha(a)=a^ib^j \qquad\textup{and}\qquad \alpha(b)=b^k \]
for some $i,j,k\in\Z$. Then $\gcd(k,r)=1$ since $\alpha(b)$ has order $r$, 
and $\gcd(i,q)=1$ since $\gen{\alpha(a),\alpha(b)}=G$. Since 
$\alpha$ is a homomorphism, 
	\[ a^{it^k}b^j = b^k(a^ib^j)b^{-k} = 
	\alpha(b)\alpha(a)\alpha(b)^{-1} = \alpha(bab^{-1}) 
	= \alpha(a^t) = (a^ib^j)^t = a^{iN}b^{tj}, \]
where $N=1+t^j+t^{2j}+\cdots+t^{(t-1)j}$. Since $\gcd(i,q)=1$, 
this implies that $N\equiv t^k$ (mod $qr$). Hence 
	\[ t^k \equiv N \equiv 1 + (1+jq) + (1+2jq) + \dots + (1+(t-1)jq) 
	=t+jq^2t/2 \pmod{qr} \]
(recall $q=t-1$). 
So $t^k\equiv t$ (mod $qr$) if $q$ is odd or $2r|q$, and 
$t^k\equiv t$ (mod $qr/2$) otherwise. Since $t=1+q$ and $t^k\equiv1+kq$ 
(mod $qr$), we now get that $k\equiv1$ (mod $r$) if $q$ is odd or $2r|q$, 
and $k\equiv1$ (mod $r/2$) otherwise. 

Thus if $q$ is odd or $2r|q$, then the only generators of $\gen{b}$ that 
are $\R$-$\Aut(G)$-conjugate to $b$ are $b$ and $b^{-1}$, so $\NN_G>0$ if 
$\varphi(r)\ge4$. Otherwise, $r$ is even, the $\R$-$\Aut(G)$-conjugacy 
class of $b$ in $\gen{b}$ contains the four elements $b^{\pm1}$ and 
$b^{(r/2)\pm1}$, and $\NN_G>0$ if $\varphi(r)>4$.
\end{proof}

We finish by looking at a few examples of finite simple groups $G$ where 
$\NN_G>0$. These are larger groups in most cases, but working with them has 
the advantage that the automorphism groups of simple groups are well known, 
and their outer automorphism groups are in most cases quite small. Also, the 
properties of conjugacy classes of elements of $G$ needed to determine 
$\NN_G$ are in many cases listed in the Atlas \cite{Atlas(1985)}. 

We start with the easiest case. 

\begin{Ex}\label{ex:An}
If $G$ is an alternating or symmetric group, then $\NN_G=0$.
\end{Ex}

\begin{proof} If $G\cong A_n$ or $\Sigma_n$, then two generators of the 
same cyclic subgroup of $G$ are always conjugate in $\Sigma_n$, and hence 
in $\Aut(G)$. So $\Q$-$\Aut(G)$-conjugate elements are also 
$\Aut(G)$-conjugate, and hence $\NN_G=0$. 
\end{proof}

We next look at the 26 sporadic simple groups, where we observe that the 
largest groups are not necessarily the ones for which $\NN_G$ or 
$\rk(\Wh(G))$ are largest. For example, the Whitehead group of $F_1$ (the 
monster) is finite, while $\NN_G>0$ when $G\cong F_2$ (the baby monster) or 
the Janko group $J_1$ (one of the smallest sporadic groups).

\begin{Ex}\label{ex:spor}
Among the sporadic simple groups $G$, 
\begin{enuma} 

\item $\Wh(G)$ is finite, and hence $\NN_G=0$, when $G$ is one of the five 
Matthieu simple groups, one of Conway's simple groups $\Co_n$ for 
$n=1,2,3$, or one of the groups $\HS$, $\McL$, $F_3$, or $F_1$; 

\item $\rk(\Wh(G))>0$ but $\NN_G=0$ when $G\cong J_2$, $\Suz$, or $\Fi_{22}$; 
and 

\item $\NN_G>0$ when $G$ is one of the Janko groups $J_n$ for $n=1,3,4$, or 
$G\cong \He$, $\Ly$, $\Ru$, $\ON$, $\Fi_{23}$, $\Fi_{24}'$, $F_5$ or $F_2$. 

\end{enuma}
\end{Ex}

\begin{proof} From the character tables in~\cite{Atlas(1985)}, we see that 
every pair of $\Q$-$G$-conjugate elements is $\R$-$G$-conjugate, and hence 
$\Wh(G)$ is finite, whenever $G$ is one of the Matthieu groups, one of 
Conway's simple groups $\Co_n$ for $n=1,2,3$, or $\HS$, $\McL$, $F_3$ or 
$F_1$. Note that an entry ``\texttt{5A B*}'' at the top of the character table 
means that class $5B$ is $\Q$-$G$-conjugate but not 
$\R$-$G$-conjugate to class $5A$, while ``\texttt{5A B**}'' means 
that the classes are $\R$-$G$-conjugate.

The computations of $\rk(\Wh(G))$ and $\NN_G$ in the other cases are 
described in Table~\ref{tbl:spor}. The actions of $\Out(G)$ on the 
conjugacy classes in $G$ can also be read off from the character tables in 
\cite{Atlas(1985)}. For example, when $G=\He$ (Held's group), the 
characters of the two classes $17A$ and $17B$ differ only for the 
characters $\chi_7$ and $\chi_8$, where $\chi_7(g)=\chi_8(h)$ when $g$ is 
in $17A$ and $h$ in $17B$ or vice versa. By the table, these two characters 
are exchanged (``fused'') by $\Out(G)\cong C_2$, and hence $\Out(G)$ must 
also exchange the classes $17A$ and $17B$. In contrast, for the same group 
$G$, the classes $21A$ and $21B$ differ only for the characters $\chi_{30}$ 
and $\chi_{31}$, these are not fused in $\Aut(G)$, and so $21A$ and $21B$ 
are not exchanged by $\Out(G)$. 

The only case where such arguments become more complicated is that of the 
four classes $16ABCD$ when $G=\ON$. Based on the action of $\Out(G)\cong 
C_2$ on the eight characters whose values are not constant on these four 
classes, one shows that the action of $\Out(G)\cong C_2$ exchanges the 
classes $16A$ and $16C$, and also the classes $16B$ and $16D$. 
\end{proof}

\begin{table}[ht] 
\begin{center}
\renewcommand{\arraystretch}{1.2}
\newcommand{\Sm}[1]{\text{\Small{$#1$}}}
\[ \begin{array}{c|ccl}
G & \rk(\Wh(G)) & \NN_G & \textup{classes} \\\hline
J_1 & 5 & 5 & 5AB,10AB,15AB,19ABC \\
J_2 & 5 & 0 & 5AB!,5CD!,10AB!,10CD!,15AB! \\
J_3 & 6 & 3 & 5AB!,9ABC,10AB!,15AB!,17AB \\
J_4 & 11 & 11 & 20AB,24AB,31ABC,33AB,37ABC,40AB,43ABC,66AB \\
\Suz & 3 & 0 & 13AB!,15AB!,21AB! \\
\He & 2 & 1 & 17AB!,21AB \\
\Ly & 11 & 11 & 21AB,24BC,31ABCDE,37AB,40AB,42AB,67ABC \\
\Ru & 7 & 7 & 14ABC,20BC,24AB,26ABC,29AB \\
\ON & 6 & 5 & 15AB,16AB,16CD,19ABC,28AB \\
\Fi_{22} & 1 & 0 & 13AB! \\
\Fi_{23} & 3 & 3 & 13AB,26AB,39AB \\
\Fi_{24}' & 8 & 2 & 21CD!,24FG!,29AB!,33AB,39AB,39CD!,42BC!,45AB! \\
F_5 & 8 & 1 & 5CD!,10DE!,10GH!,15BC!,20AB,20DE!,25AB!,30BC! \\
F_2 & 3 & 3 & 32AB,34BC,56AB \\
\end{array} \] 
\end{center}%
\caption{In the last column, we list families of $G$-conjugacy classes 
that are $\Q$-$G$-conjugate but not $\R$-$G$-conjugate to each other. For 
example, ``$19ABC$'' means that the three classes $19A$, $19B$, and $19C$ 
together form a $\Q$-$G$-conjugacy class, no two of which are 
$\R$-$G$-conjugate. An exclamation point ``!'' means that the classes are 
permuted transitively by $\Aut(G)$. When $G=\ON$, the classes $16AB$ are 
exchanged with $16CD$ by an outer automorphism.} 
\label{tbl:spor}
\end{table}

We finish by looking at the projective special linear groups.

\begin{Ex}\label{ex:Ln(q)}
Assume $G\cong\PSL_n(q)$, where $n\ge2$, $q$ is a prime power, and $G$ is 
simple. Then $\NN_G>0$, except when $G$ is one of the groups $\PSL_2(4)\cong 
\PSL_2(5)\cong A_5$, $\PSL_2(7)\cong \PSL_3(2)$, $\PSL_2(8)$, $\PSL_2(9)\cong A_6$, 
$\PSL_3(4)$, or $\PSL_4(2)\cong A_8$. 
\end{Ex}

\begin{proof} Choose an ordered $\F_q$-basis for $\F_{q^n}$, and let 
$\chi\:\F_{q^n}^\times\too\GL_n(q)$ be the injective homomorphism that 
sends an element $u\in\F_{q^n}^\times$ to the matrix for multiplication by 
$u$ on $\F_{q^n}$. Fix a generator $u_0\in\F_{q^n}^\times$, set 
$\til{x}=\chi(u_0^{q-1})\in\SL_n(q)$, and let $x\in G=\PSL_n(q)$ be its 
class modulo the center. Then $\til{x}$ has order $(q^n-1)/(q-1)$, and so 
$x$ has order 
	\[ M \defeq \frac{q^n-1}{(q-1)\cdot\gcd(q-1,n)}. \] 
Note that $\frac{q^n-1}{q-1}\equiv n$ (mod $q-1$), and hence 
$\gcd(q-1,n)=\gcd(q-1,\frac{q^n-1}{q-1})$. Also, $\Aut_G(\gen{x})\cong C_n$, 
generated by the Frobenius automorphism $(x\mapsto x^q)$. 

Assume $q=p^k$ where $p$ is prime and $k\ge1$. The 
automorphism group $\Aut(G)$ is generated by $\Inn(G)$, the field 
automorphism that sends $x$ to $x^p$, diagonal automorphisms induced by 
conjugation by $\chi(u_0)$ that send $x$ to itself, and if $n\ge3$, the 
graph automorphism ``transpose inverse'' that (up to inner automorphism) 
sends $x$ to $x^{-1}$ (see, e.g., \cite[\S\,3.3.4]{Wilson(2009)}). Thus 
$\Aut_{\Aut(G)}(\gen{x})$ is generated by the automorphisms $(x\mapsto 
x^p)$ of order $nk$ and (if $n\ge3$) $(x\mapsto x^{-1})$. It follows that
	\beqq 
	\left.\begin{array}{r} 
	\textup{$n=2$ and $\varphi(M)>2k$} \\[3pt]
	\textup{or $n\ge3$ and $\varphi(M)>2nk$}
	\end{array} \right\} 
	\implies \NN_G>0. \label{e:Lnq-1} \eeqq

\boldd{Assume $n=2$ and $p=2$.} Then $M=2^k+1$ is odd, and 
$\Aut_{\Aut(G)}(\gen{x})$ is cyclic of order $2k$, generated by the 
Frobenius automorphism $(x\mapsto x^2)$. If $M$ is divisible by two or more 
odd primes, then $\Aut(\gen{x})$ is not cyclic, and hence $\NN_G>0$. 
Otherwise, $\varphi(M)\ge\frac23\cdot M=\frac23(2^k+1)$, and so 
$\NN_G>0$ if $2^k+1>3k$. This holds whenever $k\ge4$, and so $\NN_G>0$ 
whenever $q=2^k\ge16$. 

\boldd{Assume $n=2$ and $p$ is odd.} Then $M=(q+1)/2$, and again, 
$\Aut_{\Aut(G)}(\gen{x})$ is cyclic of order $2k$. If $M$ is divisible by 
two or more odd primes, then $\NN_G>0$ since $\Aut(\gen{x})$ is not cyclic. 
Otherwise, $\varphi(M)\ge\frac12\cdot\frac23\cdot M=(q+1)/6$, and so by 
\eqref{e:Lnq-1}, $\NN_G>0$ if $q+1>12k$. This holds for $q\ge13$ when 
$k=1$, for $q=p^2\ge25$ when $k=2$, for $q=p^3\ge5^3$ when $k=3$, and for 
all odd primes $p$ when $k\ge4$. Thus by \eqref{e:Lnq-1}, $\NN_G>0$ 
whenever $q\ge13$ and $q\ne27$. 

\boldd{Assume $n\ge3$.} In these cases, $\Aut_{\Aut(G)}(\gen{x})$ is the 
product of cyclic groups of order $nk$ and $2$. If $M$ is 
divisible by three or more odd primes, then $\Aut(\gen{x})$ has $2$-rank at 
least $3$, hence is not equal to $\Aut_{\Aut(G)}(\gen{x})$, and so 
$\NN_G>0$. If $M$ is divisible by at most two odd primes, then 
$\varphi(M)\ge\frac12\cdot\frac23\cdot\frac45\cdot M=\frac4{15}M$. So by 
\eqref{e:Lnq-1}, $\NN_G>0$ if $2kn<\frac4{15}M$, and since 
$M\ge\frac1n\cdot\frac{q^n-1}{q-1}>\frac1nq^{n-1}$, we have 
	\beqq  q^{n-1} = p^{k(n-1)} 
	\ge \tfrac{15}2\cdot kn^2 
	\quad\implies\quad \NN_G>0.\label{e:Lnq-2} \eeqq

By straightforward computation, the inequality in \eqref{e:Lnq-2} holds 
whenever $(p,k,n)$ is equal to one of the triples 
\begin{small} \def\8{\hskip3mm}
	\[ (11,1,3) \8(5,2,3) \8(3,3,3) \8(2,5,3) \8(5,1,4) \8
	(3,2,4) \8(2,3,4) \8 (2,2,6) \8(3,1,7) \8(2,1,11) \]
\end{small}%
Also, if this inequality holds for a given triple $(p_0,k_0,n_0)$ (with 
$n_0\ge3$), then it holds for $(p,k,n)$ (and hence $\NN_G>0$) whenever 
$p\ge p_0$, $k\ge k_0$, and $n\ge n_0$. Among the pairs $(p^k,n)$ with 
$n\ge3$ for which \eqref{e:Lnq-2} does not hold, the inequality in 
\eqref{e:Lnq-1} holds (so $\NN_G>0$) in the following cases:
\[ \renewcommand{\arraystretch}{1.2}
\begin{array}{c|cccccccccccccccc}
q & 3 & 5 & 7 & 8 & 9 & 16 & 4 & 2 & 3 & 4 & 2 & 3 & 2 & 2 & 2 & 2 \\
n & 3 & 3 & 3 & 3 & 3 & 3 & 4 & 5 & 5 & 5 & 6 & 6 & 7 & 8 & 9 & 10 \\\hline
M & 13 & 31 & 19 & 73 & 91 & 91 & 85 & 31 & 121 & 341 & 63 & 182 & 127 & 255 & 511 
& 1023 \\
\varphi(M) & 12 & 30 & 18 & 72 & 72 & 72 & 64 & 30 & 110 & 300 & 36 & 72 & 126 & 
128 & 432 & 600 \\
2kn & 6 & 6 & 6 & 18 & 12 & 24 & 16 & 10 & 10 & 20 & 12 & 12 & 14 & 16 & 18 & 20
\end{array} \]

\boldd{Remaining cases (all $n$).} We are left with the groups $\PSL_2(q)$ for 
$q\le11$ or $q=27$, $\PSL_3(q)$ for $q=2,4$, and $\PSL_4(q)$ for $q=2,3$. 
The groups $\PSL_2(4)\cong \PSL_2(5)\cong A_5$, $\PSL_2(9)\cong A_6$, and 
$\PSL_4(2)\cong A_8$ were handled in Example~\ref{ex:An}. Also, 
$\PSL_3(2)\cong \PSL_2(7)$, so it remains to look at $\PSL_2(7)$, 
$\PSL_2(8)$, $\PSL_2(11)$, $\PSL_2(27)$, $\PSL_3(4)$, and $\PSL_4(3)$. In 
each of these cases, the result follows with the help of the character 
tables in~\cite{Atlas(1985)}.

For example, when $G=\PSL_2(11)$, the classes $5A$ and $5B$ are not 
$\R$-$\Aut(G)$-conjugate, so $\NN_G>0$. When $G=\PSL_2(27)$, there are six 
conjugacy classes of elements of order $13$ (all of them 
$\Q$-$G$-conjugate), permuted in two orbits by 
field automorphisms of order $3$, hence forming two 
$\R$-$\Aut(G)$-conjugacy classes. When $G=\PSL_4(3)$, the subgroups of 
order $13$ have automizers in $\Aut(G)$ of order $6$, so their generators 
form two $\R$-$\Aut(G)$-conjugacy classes. Thus $\NN_G>0$ in the last two 
cases. When $G=\PSL_2(7)$, $\PSL_2(8)$, or $\PSL_3(4)$, each 
$\Q$-$\Aut(G)$-conjugacy class is permuted transitively by $\Aut(G)$, and 
so $\NN_G=0$. 
\end{proof}

\bigskip

Conflict of interest statement: On behalf of all authors, the corresponding author states that there is no conflict of interest.

\bigskip

Data availability statement: No datasets were generated or analysed during the current study.



\end{document}